\documentclass[reqno]{article}
\usepackage{amsfonts, amsmath}
\begin{document}
\setlength{\textwidth}{16cm}
\renewcommand{\thesection}{\arabic{section}.}
\renewcommand{\theequation}{\arabic{section}.\arabic{equation}}
\newcommand{\be}{\begin{eqnarray}}
\newcommand{\en}{\end{eqnarray}}
\newcommand{\no}{\nonumber}
\newcommand{\la}{\lambda}
\newcommand{\ep}{\epsilon}
\newcommand{\lan}{\langle}
\newcommand{\ra}{\rangle}
\newcommand{\de}{\delta}
\newcommand{\ov}{\overline}
\newcommand{\bet}{\beta}
\newcommand{\al}{\alpha}
\newcommand{\p}{\partial}
\newcommand{\fr}{\frac}
\newcommand{\we}{\wedge}
\newcommand{\D}{\Delta}
\newcommand{\R}{\mathbb{R}}
\newcommand{\ri}{\rightarrow}
\newcommand{\m}{\mathbb}
\newcommand{\om}{\Omega}
\newcommand{\na}{\nabla}
\newcommand{\vs}{\vskip0.3cm}
\newcommand{\pa}{\partial}
\newcommand{\va}{\varsigma}
\pagestyle{myheadings}
\newtheorem{theorem}{Theorem}[section]
\newtheorem{lemma}[theorem]{Lemma}
\newtheorem{condition}[theorem]{Conditions}
\newtheorem{corollary}[theorem]{Corollary}
\newtheorem{proposition}[theorem]{Proposition}
\newtheorem{remark}[theorem]{Remark}
\newtheorem{example}[theorem]{Example}
\newtheorem{conjecture}[theorem]{Conjecture}
\newtheorem{hypotheses}[theorem]{Hypotheses}
\newtheorem{defn}[theorem]{Definition}
\newtheorem{maintheorem}{Main Theorem}[section]
\renewcommand{\thefootnote}{}

\title{On Ashbaugh-Benguria's Conjecture about Lower Order Dirichlet  Eigenvalues of the Laplacian}
\footnotetext{ MSC 2010: 35P15; 58C40.
\\
\ \ \ \ \ \hspace*{4ex} Key Words: Ashbaugh-Benguria's Conjecture,
Isoperimetric Inequality, Eigenvalues, Dirichlet Problem.}
\author{
Qiaoling Wang$^{a}$, Changyu Xia$^{a}$}

\date{}

\maketitle ~~~\\[-15mm]

\begin{center}
{\footnotesize  $a$. Departamento de Matem\'atica, Universidade de
Brasilia, 70910-900-Brasilia-DF, Brazil Email: wang@mat.unb.br(Q.
Wang), xia@mat.unb.br(C. Xia).  }
\end{center}

\begin{abstract}
We prove an isoperimetric inequality for lower order eigenvalues of
the Dirichlet Laplacian on bounded domains of a Euclidean space
which strengthens the celebrated Ashbaugh-Benguria inequality
conjectured by Payne-P\'olya-Weinberger on the ratio of the first
two Dirichlet eigenvalues  and makes an important step toward the
proof of a conjecture by Ashbaugh-Benguria.
 \end{abstract}

\markright{\sl\hfill  Q. Wang, C. Xia \hfill}

\section{Introduction}
\renewcommand{\thesection}{\arabic{section}}
\renewcommand{\theequation}{\thesection.\arabic{equation}}
\setcounter{equation}{0}

Let $\om $ be a bounded domain in $\R^n$,  $n\geq 2$. Let us denote
by $\D$ the Laplace operator on $\R^n$ and consider the homogeneous
membrane problem
 \be\label{1.1}\left\{\begin{array}{l}
\Delta u = -\la u \ \ \ {\rm in  \ }\ \ \om,\\
\ \ \ u =  0 \ \ \ \ \ \ \ {\rm on  \ }\ \pa\om.
\end{array}\right.
\en It is well known that the spectrum of  (1.1) is real and
discrete consisting in a sequence \be \no 0<\la_1<
\la_2\leq\la_3\cdots \rightarrow +\infty,\en where each eigenvalue
is repeated with its multiplicity. An important issue in analysis
and geometry is to give good estimates to these and other
eigenvalues, especially to obtain isoperimetric bounds for them.
When $\om=\mathbb{B}^n$ is the $n$-dimensional unit ball in $\R^n$,
it is well known that $\la_1(\m{B}^n)=j^2_{n/2-1,1}$ and
$\la_2(\m{B}^n)=\cdots=\la_{n+1}(\m{B}^n)=j^2_{n/2,1}$, where
$j_{p,k}$ denotes the $k$th positive zero of the Bessel function
$J_p(x)$ of the first kind of order $p$. One of the earliest
isoperimetric inequalities for an eigenvalue is  the Faber-Krahn
inequality \cite{f,k1,k2}  conjectured by Rayleigh \cite{r} in 1877:
\be \la_1(\om)\geq\left(\fr{|\m{B}^n|}{|\om|}\right)^{\fr
2n}j^2_{\fr n2-1,1}, \en with equality if and only if $\om$
 is an $n$-ball. Here, $|\om|$ denotes the volume of $\om$.
In 1956, Payne-P\'olya-Weinberger proposed the following well-known
conjecture \cite{ppw1}: \vskip0.3cm {\bf Payne-P\'olya-Weinberger
Conjecture.} {\it The eigenvalues of (1.1) satisfy \be
\fr{\la_2(\om)}{\la_1(\om)}\leq \fr{\la_2(\m{B}^n)}{\la_1(\m{B}^n)},
\en \be \fr{\la_2(\om)+\cdots+\la_{n+1}(\om)}{\la_1(\om)}\leq
n\fr{\la_2(\m{B}^n)}{\la_1(\m{B}^n)}. \en }

The conjecture (1.3) was studied by many mathematicians, for
examples, Payne, P\'olya and Weinberger \cite{ppw1,ppw2}, Brands
\cite{b},  Chiti \cite{ch1,ch2}, de Vries \cite{d}, Hile and Protter
\cite{hp}. Finally, Ashbaugh and Benguria proved this conjecture
\cite{ab1,ab2,ab3}.  Ashbaugh-Benguria \cite{ab6} and Benguria-Linde
\cite{bl} also proved similar inequalities for the first $(n+1)$
Dirichlet eigenvalues of the Laplacian on bounded domains in a
hemisphere or a hyperbolic space.

The conjecture (1.4) is stronger than (1.3) and was also studied by
many authors. In 1956, Payne, P\'olya and Weinberger \cite{ppw2}
proved that for $\om\subset\R^2$, \be \fr{\la_2+\la_3}{\la_1}\leq 6,
\en which was improved by Brands \cite{b} to \be
\fr{\la_2+\la_3}{\la_1}\leq 3+\sqrt 7. \en Furthermore, Hile-Protter
\cite{hp} obtained \be \fr{\la_2+\la_3}{\la_1}\leq 5.622.\en  In
\cite{m}, Marcellini obtained the bound \be
\fr{\la_2+\la_3}{\la_1}\leq\fr{15+\sqrt{345}}6.\en
 Chen-Zheng proved
in \cite{cz} \be \fr{\la_2+\la_3}{\la_1}\leq 5.3507^-. \en

For general dimensions $n\geq 2$, Thompson \cite{t} obtained the
bound (see also \cite{ab5})
 \be \fr{\la_1+\la_2\cdots
\la_{n+1}}{\la_1}\leq (n+4).\en In \cite{ab5}, Ashbaugh-Benguria
proved \be \fr 1{\la_2-\la_1}+\cdots +\fr
1{\la_{n+1}-\la_1}\geq\fr{2j^2_{\fr n2-1, 1}+n(n-4)}{6\la_1}. \en
They observed that \cite{ab5} \be \fr{2j^2_{\fr n2-1,
1}+n(n-4)}6\sim\fr{n^2}4\left[1+\fr 23(1.8557571)\fr{2^{\fr
23}}{n^{\fr 23}}-\fr 4n + O(n^{-\fr 43})\right], \en whereas \be \fr
n{\left(\fr{j_{\fr n2}}{j_{\fr
n2-1}}\right)^2-1}\sim\fr{n^2}4\left[1+\fr 23(1.8557571)\fr{2^{\fr
23}}{n^{\fr 23}}-\fr 2n + O(n^{-\fr 43})\right] \en and also
conjectured that \cite{ab4,ab5}\be
\fr{\la_1}{\la_2-\la_1}+\cdots+\fr{\la_1}{\la_{n+1}-\la_1}\geq \fr
n{\left(\fr{j_{\fr n2}}{j_{\fr n2-1}}\right)^2-1}\en with equality
if and only if $\om$
 is an $n$-ball.

\vskip0.3cm Ashbaugh \cite{a} and Henrot \cite{h} mentioned this
conjecture again. One can also formulate a similar conjecture for
the first $(n+1)$ eigenvalues of the Dirichlet Laplacian on bounded
domains in a hemisphere or a hyperbolic space.

In this paper, we prove the following isoperimetric inequality which
supports strongly the  conjecture (1.14).

\begin{theorem}\label{th1} Let $\om$ be a bounded domain  with smooth boundary in
$\R^n$. Then the first $n$ Dirichlet eigenvalues of $\om$ satisfy
\be\fr{\la_1}{\la_2-\la_1}+\cdots+\fr{\la_1}{\la_{n}-\la_1}\geq \fr
{n-1}{\left(\fr{j_{\fr n2}}{j_{\fr n2-1}}\right)^2-1},\en with
equality holding if and only if $\om$
 is an $n$-ball.
\end{theorem}

For eigenvalues $0=\mu_0<\mu_1\leq\mu_2\leq\cdots \ri +\infty$ of
the Neumann problem
\begin{eqnarray}\label{a1}
\left\{\begin{array}{ccc} \D u=\mu u\,&&~\mbox{in} ~~ \om, \\[2mm]
\frac{\p u}{\p\nu}=0, &&~~\mbox{on}~~\partial \om,
\end{array}\right.
\end{eqnarray}
 where $\frac{\p}{\p\nu}$ is the outer normal derivative, the well-known Szeg\"o-Weinberger inequality
 states that \cite{s,w}
  \be \label{int2} \mu_1(\om)|\om|^{2/n}\leq
\mu_1(\mathbb{B}^n)|\mathbb{B}^n|^{2/n}, \en with equality holding
if and only if $\om$ is a ball in $\mathbb{R}^n$.
  Ashbaugh and Benguria  conjectured in \cite{ab4}
 that
\begin{eqnarray}\label{a2}
\sum_{i=1}^{n}\frac{1}{\mu_i(\om)}\geq\frac{n}{\mu_1(B_\Omega)},
~~\mathrm{with~equality~if~and~only~if}~\Omega~\mathrm{is~a~ball},
\end{eqnarray}
where $B_{\om}\subset \R^n$ is a ball of same volume as $\om$. In
\cite{wx}, the authors proved the following  inequality
\begin{eqnarray}\label{a3}
\sum_{i=1}^{n-1}\frac{1}{\mu_i(\om)}\geq\frac{n-1}{\mu_1(B_\Omega)},
~~\mathrm{with~equality~if~and~only~if}~\Omega~\mathrm{is~a~ball},
\end{eqnarray} which supports this
conjecture of  Ashbaugh and Benguria.

\section{A proof of Theorem 1.1.}
\setcounter{equation}{0} Before proving Theorem \ref{th1}, let us
recall some known facts (Cf. \cite{ab1,ab2,ab3,c,h,sy,x}). Let
$\{u_j\}_{j=1}^{\infty}$ be an orthonormal set of eigenfunctions of
the problem (1.1), that is, \be\label{pth1.6}
\left\{\begin{array}{l} \Delta u_i= -\la_i u_i \ \ \ {\rm in} \ \ \
\om,\\
 u_i|_{\pa\om}=0,\\
\int_{\om} u_i u_j dx=\delta_{ij}.
\end{array}\right.,  \en where
$dx$ denotes the volume element of $\om$. For each $k=1,2,\cdots,$
the variational characterization of $\la_{k+1}(\om)$ is given by
 \be\label{pth0.5}
\la_{k+1}(\om)=\underset{\underset{\int_{\om}\phi u_i dx=0,
i=1,\cdots,k}{\phi\in
H_0^1(\om)\setminus\{0\}}}{\inf}\fr{\int_{\om}|\na\phi|^2
dx}{\int_{\om} \phi^2 dx}.
 \en
Let $B_r$ be a ball of radius $r$ centered at the origin in $\R^n$.
It is known that \be \la_1(B_r)=\left(\fr{j_{(n-2)/2,1}}r\right)^2
\en with its corresponding eigenfunction given by the radial
function \be u(x):=c|x|^{1-\fr n2}J_{\fr
n2-1}\left(\fr{j_{(n-2)/2,1}}r|x|\right), \en where $c$ is a nonzero
constant. The second Dirichlet eigenvalue of $B_r$
 has multiplicity $n$, that is,
 \be
 \la_2(B_r)=\cdots
=\la_{n+1}(B_r)=\fr{j_{n/2,1}^2}{r^2} \en and a basis for the
eigenspace corresponding to $\la_2(B_r)$ consists of
\be\label{pth1.2} \xi_i(x)=|x|^{1-\fr n2} J_{n/2}\left(\fr{j_{n/2,
1}|x|}r\right)\fr{x_i}{|x|}, \ \ i=1,\cdots, n. \en Define a
function $w: [0, +\infty)\ri \R$ by
\be\label{1.1}w(t)\equiv\left\{\begin{array}{l}
\fr{J_{\fr n2}(\beta t)}{J_{\fr n2-1}(\alpha t)}  \ \ \ \ \ \ \ \ \ \ \ \ \ \ \ \ \ \ \ \ \ \ \ {\rm for  \ }\ 0\leq t<1, \\
\ \ \ w(1) \equiv \lim_{t\ri 1^-} w(t) \ \ \ \ \  {\rm for  \ }\
t\geq 1,
\end{array}\right.
\en where $\alpha=j_{n/2-1,1}, \beta=j_{n/2,1}.$ We have $w(0)=0,
w(t)>0, \ \forall t\in (0, +\infty)$  and for any $t\geq 0$, one
concludes from  Theorem 3.3 in \cite{ab2} that
 \be (w^{\prime}(t))^2\leq
\left(\fr{w(t)}{t}\right)^2. \en Let $\gamma =\sqrt{\la_1}/\alpha$
and set \be B(t)\equiv w^{\prime}(t)^2+(n-1)\fr{w(t)^2}{t^2}; \en
then( Cf. (2.14), (2.15) and (2.22) in \cite{ab2}) \be
\fr{\int_{\om} B(\gamma |x|)u_1^2 dx}{\int_{\om} w(\gamma |x|)^2
u_1^2 dx}\leq\beta^2- \alpha^2. \en

{\it Proof of Theorem \ref{th1}.} Observe that if \be Q\neq 0 \ {\rm
and\ } \int_{\om} Qu_1^2 dx=\int_{\om}Qu_1u_2dx=\cdots =\int_{\om}
Qu_1 u_k dx=0, \en then (2.2) gives
 \be \la_{k+1}\leq \fr{\int_{\om}|\na(Qu_1)|^2  dx}{\int_{\om} Q^2 u_1^2 dx}, \en
which, yields by integration by parts that
 \be \la_{k+1}-\la_1\leq
\fr{\int_{\om}|\na Q|^2 u_1^2 dx}{\int_{\om} Q^2 u_1^2 dx}. \en
 We define $g: [0, +\infty)\ri\R$ by
\be\label{pth1.6} g(t)= w(\gamma t) \en and fix an orthonormal basis
$\{e_i\}_{i=1}^n$ of $\R^n$. By using the Brouwer fixed-point
theorem, we can  choose the origin of $\R^n$ so that (Cf.
\cite{ab2}) \be \label{pth1.7} \int_{\om}\langle x,
e_i\rangle\fr{g(|x|)}{|x|} u_1^2 dx =0,\ \ i=1,\cdots, n.\en Next we
show that there exists a new orthonormal basis
$\{e_i^{\prime}\}_{i=1}^n$ of $\R^n$ such that \be\label{pth1.8}
\int_{\om}\langle x, e_i^{\prime}\rangle\fr{g(|x|)}{|x|}u_1 u_{j+1}
dx=0, \en for $j=1,\cdots, i-1$ and $i=2,\cdots, n$. To see this, we
define an $n \times n$ matrix $P=\left(p_{ij}\right)$ by \be
p_{ij}=\int_{\om} \langle x, e_i\rangle\fr{g(|x|)}{|x|}u_1 u_{j+1}
dx, \ i,j=1,2,\cdots,n.\en Using the orthogonalization of Gram and
Schmidt (QR-factorization theorem), one can find an upper triangle
matrix $T=(T_{ij})$ and an orthogonal matrix $U=(a_{ij})$ such that
$T=UP$. Hence,
\begin{eqnarray*}
T_{ij}=\sum_{k=1}^n a_{ik}p_{kj}=\int_{\om} \sum_{k=1}^n
a_{ik}\langle x, e_k\rangle\fr{g(|x|)}{|x|}u_1 u_{j+1} dx =0,\ \
1\leq j<i\leq n.
\end{eqnarray*}
Letting $e_i^{\prime}=\sum_{k=1}^n  a_{ik}e_k, \ i=1,...,n$, one
gets  (\ref{pth1.8}). Let us denote by $x_1, x_2,\cdots, x_n$ the
coordinate functions of $\R^n$ with respect to the base
$\{e_i^{\prime}\}_{i=1}^n$, that is, $x_i=\langle x,
e_i^{\prime}\rangle, \ x\in\R^n$. From (\ref{pth1.7}) and
(\ref{pth1.8}), we have \be\label{pth1.9} \int_{\om}
g(|x|)\fr{x_i}{|x|}u_1 u_{j+1} dx=0, \ i=1,\cdots, n, \ j=0,\cdots,
i-1. \en Let
$$\phi_k=g(|x|)\fr{x_k}{|x|},\ \  k=1,\cdots, n;$$
 then \be
\phi_k\not\equiv 0 {\ \ \rm and\ \ } \int_{\om}\phi_k u_1^2
dx=\cdots=\int_{\om}\phi_k u_1 u_k dx=0. \en It then follows from
(2.13) that
\begin{eqnarray}\label{pth1.11}
(\la_{k+1}-\la_1)\int_{\om} \phi_k^2 u_1^2 dx \leq\int_{\om}
|\na\phi_k|^2 u_1^2 dx,\ k=1,\cdots,n. \en Substituting \be\label{p}
|\na \phi_k|^2&=& g^{\prime}(|x|)^2\fr{x_k^2}{|x|^2} +
\fr{g(|x|)^2}{|x|^2}\left(1-\fr{x_k^2}{|x|^2}\right)\\ \no &=&
\left(
g^{\prime}(|x|)^2-\fr{g(|x|)^2}{|x|^2}\right)\fr{x_k^2}{|x|^2}+\fr{g(|x|)^2}{|x|^2}\en
into (2.20) and dividing by $(\la_{k+1}-\la_1),$ we have for
$k=1,\cdots,n,$ that \be\label{pth1.12} \int_{\om}\phi_k^2 u_1^2
dx&\leq&\fr 1{\la_{k+1}-\la_1}\int_{\om}\left(
g^{\prime}(|x|)^2-\fr{g(|x|)^2}{|x|^2}\right)\fr{x_k^2}{|x|^2}u_1^2dx\\
\no \ \ & & + \fr
1{\la_{k+1}-\la_1}\int_{\om}\fr{g(|x|)^2}{|x|^2}u_1^2dx. \en Summing
on $k$ from $1$ to $n$, one gets \be\int_{\om} g(|x|)^2u_1^2
dx&\leq& \sum_{k=1}^n \fr
1{\la_{k+1}-\la_1}\int_{\om}\fr{g(|x|)^2}{|x|^2}u_1^2dx\\
\no & & +\sum_{k=1}^n\fr 1{\la_{k+1}-\la_1}\int_{\om}\left(
g^{\prime}(|x|)^2-\fr{g(|x|)^2}{|x|^2}\right)\fr{x_k^2}{|x|^2}u_1^2dx.\en
Observe
that \be& &  \sum_{k=1}^n \fr 1{\la_{k+1}-\la_1}\fr{x_k^2}{|x|^2}\\
\no &=& \sum_{k=1}^{n-1} \fr 1{\la_{k+1}-\la_1}\fr{x_k^2}{|x|^2}+\fr
1{\la_{n+1}-\la_1}\fr{x_n^2}{|x|^2}\\ \no &=&
 \sum_{k=1}^{n-1} \fr 1{\la_{k+1}-\la_1}\fr{x_k^2}{|x|^2}+\fr
1{\la_{n+1}-\la_1}\left(1-\sum_{k=1}^{n-1}\fr{x_k^2}{|x|^2}\right).
\en Therefore, \be& & \sum_{k=1}^n\fr
1{\la_{k+1}-\la_1}\int_{\om}\left(
g^{\prime}(|x|)^2-\fr{g(|x|)^2}{|x|^2}\right)\fr{x_k^2}{|x|^2}u_1^2 dx\\
\no &=& \sum_{k=1}^{n-1}\fr 1{\la_{k+1}-\la_1}\int_{\om}\left(
g^{\prime}(|x|)^2-\fr{g(|x|)^2}{|x|^2}\right)\fr{x_k^2}{|x|^2}u_1^2 dx\\
\no & & +\fr 1{\la_{n+1}-\la_1}\int_{\om}\left(
g^{\prime}(|x|)^2-\fr{g(|x|)^2}{|x|^2}\right)u_1^2 dx\\
\no & & -\fr 1{\la_{n+1}-\la_1}\int_{\om}\left(
g^{\prime}(|x|)^2-\fr{g(|x|)^2}{|x|^2}\right)\sum_{k=1}^{n-1}\fr{x_k^2}{|x|^2}
u_1^2 dx.\\ \no &=& \sum_{k=1}^{n-1}\int_{\om}\left(\fr
1{\la_{k+1}-\la_1}-\fr 1{\la_{n+1}-\la_1}\right)\left(
g^{\prime}(|x|)^2-\fr{g(|x|)^2}{|x|^2}\right)\fr{x_k^2}{|x|^2}u_1^2
dx\\ \no & &  +\fr 1{\la_{n+1}-\la_1}\int_{\om}\left(
g^{\prime}(|x|)^2-\fr{g(|x|)^2}{|x|^2}\right)u_1^2 dx. \en We have
\be\no \fr 1{\la_{k+1}-\la_1}-\fr 1{\la_{n+1}-\la_1}\geq 0,\
k=1,\cdots,n-1. \en It follow from (2.8) and (2.14) that \be\no
g^{\prime}(|x|)^2-\fr{g(|x|)^2}{|x|^2}\leq 0 \ \ {\rm on} \ \om. \en
Thus, \be\no \sum_{k=1}^{n-1}\int_{\om}\left(\fr
1{\la_{k+1}-\la_1}-\fr 1{\la_{n+1}-\la_1}\right)\left(
g^{\prime}(|x|)^2-\fr{g(|x|)^2}{|x|^2}\right)\fr{x_k^2}{|x|^2}u_1^2
dx\leq 0, \en which, combining with (2.23) and (2.25), gives
 \be\int_{\om} g(|x|)^2 dx&\leq& \fr
1{\la_{n+1}-\la_1}\int_{\om}\left(
g^{\prime}(|x|)^2-\fr{g(|x|)^2}{|x|^2}\right)u_1^2dx \\
\no & & \ \  + \sum_{k=1}^n \fr
1{\la_{k+1}-\la_1}\int_{\om}\fr{g(|x|)^2}{|x|^2}u_1^2 dx\\ \no &=&
 \fr 1{\la_{n+1}-\la_1}\int_{\om} g^{\prime}(|x|)^2u_1^2 dx +\sum_{k=1}^{n-1} \fr
1{\la_{k+1}-\la_1}\int_{\om}\fr{g(|x|)^2}{|x|^2}u_1^2 dx\\ \no
&\leq& \fr 1{n-1}\sum_{k=1}^{n-1} \fr
1{\la_{k+1}-\la_1}\int_{\om}\left(g^{\prime}(|x|)^2+(n-1)\fr{g(|x|)^2}{|x|^2}\right)u_1^2
dx. \en Consequently, we have from (2.9), (2.10), (2.14) and (2.26)
that
 \be\no
\fr 1{n-1}\sum_{k=1}^{n-1} \fr 1{\la_{k+1}-\la_1}&\geq&\fr
{\int_{\om}g(|x|)^2u_1^2
dx}{\int_{\om}\left(g^{\prime}(|x|)^2+(n-1)\fr{g(|x|)^2}{|x|^2}\right)u_1^2
dx}\\ \no &=&\fr{\alpha^2}{\la_1} \fr{\int_{\om}w(\gamma|x|)^2 u_1^2
dx}{\int_{\om} B(\gamma|x|)u_1^2 dx}\\ \no &\geq&
\fr{\alpha^2}{\la_1}\fr 1{(\beta^2-\alpha^2)}, \en which proves
(1.15). Also, one can see that the equality holds in (1.15) if and
only if $\om$ is an $n$-ball. This completes the proof of Theorem
1.1.
\section{Lower Order Dirichlet eigenvalues of general
elliptic equations} By using the arguments in the proof of Theorem
1.1 and the work of Ashbaugh-Benguria [3] one  can generalize the
inequality (1.15) to  the first $n$ eigenvalues of the following
general problem \setcounter{equation}{0} \be \left\{\begin{array}{l}
-\sum_{i,j}\fr{\partial}{\partial x_i}\left(a_{ij}(x)\fr{\partial
u}{\partial x_j}\right)+ q(x) u= \la r(x) u
\ \ {\rm in \ } \om, \\
u|_{\partial \om}=0, \end{array} \right.
 \en
where $\om$ is bounded domain with smooth boundary in $\R^n$ and
$[a_{ij}(x)]$ is symmetric positive definite for any $x\in\om$.
Namely, we have
\begin{theorem} For equation (3.1), assume that
$q\geq 0$ on $\om$ and that there are positive numbers $a$, $A$, $c$
and $C$ such that the matrix $ [a_{ij}]$ satisfies \be a\leq
[a_{ij}]\leq A \en
 in the sense of quadratic forms throughout $\om$
and \be c\leq r(x)\leq C \ \ {\rm on}\ \ \om.\en
 Then the first $n$
eigenvalues of the problem (3.1) satisfy \be
\fr{\la_1}{\la_2-\la_1}+\cdots+\fr{\la_1}{\la_{n}-\la_1}\geq \fr
{(n-1)ac}{AC\left(\left(\fr{j_{\fr n2}}{j_{\fr
n2-1}}\right)^2-1\right)}. \en Furthermore, equality holds if and
only if $c=C$, $a=A, q\equiv 0$, and $\om$ is a ball in $\R^n$.
\end{theorem}

{\it Proof.} Let  $\{v_k\}_{k=1}^{\infty}$ be an orthonormal set of
eigenfunctions of the problem (3.1), that is, \be\label{pth1.6}
\left\{\begin{array}{l} -\sum_{i,j}\fr{\partial}{\partial
x_i}\left(a_{ij}(x)\fr{\partial v_k}{\partial x_j}\right)+ q(x) v_k=
\la_k r(x) v_k
\ \ {\rm in \ } \om, \\
v_k|_{\partial \om}=0, \\
\int_{\om} r v_k v_l dx =\delta_{kl}. \end{array} \right. \en For
each $k=1,2,\cdots,$ the variational characterization of $\la_{k+1}$
of the problem (3.1) is given by
 \be\label{pth0.5}
\la_{k+1}(\om)=\underset{\underset{\int_{\om}r\psi v_i dx=0,
i=1,\cdots,k}{\psi\in
H_0^1(\om)\setminus\{0\}}}{\inf}\fr{\int_{\om}\left(\sum_{i,j}a_{ij}(x)\fr{\partial
\psi}{\partial x_i}\fr{\partial \psi}{\partial x_j}+ q\psi^2\right)
dx}{\int_{\om} r\psi^2 dx}.
 \en
Thus if $Q$ is such that $Q\neq 0$ and \be
 \int_{\om} rQ v_1^2 dx=\cdots=\int_{\om}rQv_1 v_{k}dx=0,
 \en
then
 \be
\la_{k+1}&\leq&\fr{\int_{\om}\left(\sum_{i,j}a_{ij}(x)\fr{\partial
(Qv_1)}{\partial x_i}\fr{\partial(Qv_1)}{\partial x_j} + qQ^2
v_1^2\right) dx}{\int_{\om} r(x)Q^2 v_1^2 dx}. \en It then follows
from integration by parts, (3.2), (3.3) and the fact that $v_1$ is
an eigenfunction corresponding to the eigenvalue $\la_1$ that \be
\la_{k+1}-\la_1\leq
\fr{\int_{\om}\left(\sum_{i,j}a_{ij}(x)\fr{\partial Q}{\partial
x_i}\fr{\partial Q}{\partial x_j}v_1^2\right) dx}{\int_{\om} r(x)Q^2
v_1^2 dx}\no &\leq& \fr{A}c\fr{\int_{\om}|\na Q|^2 v_1^2
dx}{\int_{\om} Q^2 v_1^2 dx}, \en Let $\omega$ and $B$ be as in
Section 2 and set \be \gamma =\fr 1{j_{n/2-1,1}}\sqrt{\fr{C\la_1}a},
\ \ g(t)=\omega(\gamma t), \ t\geq 0.\en From the proof of Theorem
1.1, we know that one can choose the origin and the coordinate
system of $\R^n$ properly so that
 \be
\int_{\om} r g(|x|)\fr{x_i}{|x|}v_1 v_{j+1} dx=0, \ i=1,\cdots, n, \
j=0,\cdots, i-1 \en and so for $k=1,\cdots,n,$
 \be
(\la_{k+1}-\la_1)\int_{\om} \left(g(|x|)\fr{x_i}{|x|}\right)^2 v_1^2
dx\leq
\fr{A}a\int_{\om}\left|\na\left(g(|x|)\fr{x_i}{|x|}\right)\right|^2
v_1^2 dx. \ \ \en Dividing (3.11) by $(\la_{k+1}-\la_1)$ and summing
on $k$ from $1$ to $n$, one gets as in Section 2 that \be
\no\int_{\om} g(|x|)^2 v_1^2 dx \leq \fr 1{n-1}\sum_{k=1}^{n-1}\fr
1{\la_{k+1}-\la_1}\fr{A}c\int_{\om}\left(g^{\prime}(|x|)^2+(n-1)\fr{g(|x|)^2}{|x|^2}\right)
v_1^2 dx. \en Hence, \be \fr 1{n-1}\sum_{k=1}^{n-1}\fr
1{\la_{k+1}-\la_1}&\geq& \fr{c}A\fr{\int_{\om} g(|x|)^2 v_1^2
dx}{\int_{\om}\left(g^{\prime}(|x|)^2+(n-1)\fr{g(|x|)^2}{|x|^2}\right)
v_1^2 dx}
\\ \no &=& \fr{c}A\fr 1{\gamma^2}\fr{\int_{\om}w(\gamma |x|)^2 v_1^2 dx}{\int_{\om}B(\gamma |x|)v_1^2
dx}\\ \no  &=& \fr{c}A\fr{a j_{n/2-1,
1}^2}{C\la_1}\fr{\int_{\om}w(\gamma |x|)^2 v_1^2
dx}{\int_{\om}B(\gamma |x|)v_1^2 dx}. \en The proof of Theorem 4.1
in \cite{ab2} shows that \be \fr{\int_{\om}w(\gamma |x|)^2 v_1^2
dx}{\int_{\om}B(\gamma |x|)v_1^2 dx}\geq \fr 1{j_{n/2,1}^2-j_{n/2-1,
1}^2}. \en Combining (3.11) and (3.12), we get (3.2). Also, we can
see that the equality  holds in (3.2) if an only if $c=C$, $a=A,
q\equiv 0$, and $\om$ is a ball in $\R^n$. This completes the proof
of Theorem 3.1.

\section*{Acknowledgments}

 Q. Wang was partially supported by CNPq,
Brazil (Grant No. 307089/2014-2). C. Xia was partially supported by
CNPq, Brazil (Grant No. 306146/2014-2).

\end{document}